\RequirePackage{fix-cm}
\documentclass[twocolumn]{svjour3}          % twocolumn
\smartqed  % flush right qed marks, e.g. at end of proof

\usepackage{url}
\usepackage{bbold}
\usepackage{graphicx}
\usepackage{amsmath}
\usepackage{tabularx}
\usepackage{mwe}
\usepackage{hyperref}
\usepackage[export]{adjustbox}
\usepackage{color,soul}
\usepackage[ruled,vlined]{algorithm2e}
\usepackage{cite}
\usepackage{float}

\begin{document}

\title{Forcing the escape: Partial control of escaping orbits from a transient chaotic region}
%\title{Partial control of the escape of a transient chaotic region}
% Force line breaks with \\

\author{Gaspar Alfaro \and 
        Rub{\'e}n Cape{\'a}ns \and
        Miguel A.F. Sanju{\'a}n
}

\institute{G.~Alfaro $\cdot$ R.~Cape{\'a}ns $\cdot$ M.A.F. Sanju{\'a}n \at
              Nonlinear Dynamics, Chaos and Complex Systems Group, Departamento de  F\'isica, Universidad Rey Juan Carlos, M\'ostoles, Madrid, Tulip\'an s/n, 28933, Spain \\
              \email{miguel.sanjuan@urjc.es}           %  \\
}

%\affiliation{Nonlinear Dynamics, Chaos and Complex Systems Group, Departamento de  F\'isica, Universidad Rey Juan Carlos, M\'ostoles, Madrid, Tulip\'an s/n, 28933, Spain}

\date{Received: date / Accepted: date   \today }

\maketitle

\begin{abstract}
A new control algorithm based on the partial control method has been developed. The general situation we are considering is an orbit starting in a certain phase space region $Q$ having a chaotic transient behavior affected by noise, so that the orbit will definitely escape from $Q$ in an unpredictable number of iterations. Thus, the goal of the algorithm is to control in a predictable manner when to escape. While partial control has been used as a way to avoid escapes, here we want to adapt it to force the escape in a controlled manner. We have introduced new tools such as escape functions and escape sets that once computed makes the control of the orbit straightforward. We have applied the new idea to three different cases in order to illustrate the various application possibilities of this new algorithm.
\keywords{Controlling chaos, partial control, transient chaos, escaping orbits}
\end{abstract}

\section{Introduction}

Even though chaotic systems are difficult to deal with due to its intrinsic unpredictability, there are nonetheless methods that allow us to control them. Different techniques for controlling chaos have been developed in the past few years. A rough classification may divide them between feedback control and non-feedback control methods. Among the first, we can consider the OGY \cite{OGY} or the Pyragas \cite {Pyragas} control methods, while on the latter random, chaotic or periodical signals are used as an appropriate mechanism to control the system.

The partial control method,  which is a feedback method, has been used in previous works \cite{zasayo08,zasa09,sazasayo12,sasayo12,capeans2014less,casasa17,casasa18}, and is applied to a map defined in a certain region $Q$ where there is transient chaos and in absence of any control the orbit will eventually escape from the region after a certain number of iterations. Furthermore, the map is subjected to a disturbance which is always larger than the applied control. The goal is to use the minimum control to keep the orbits inside $Q$ in presence of the disturbance. Precisely, one advantage of this method is the capacity to keep small the amount of control. As is well known, transient chaos is the physical manifestation of the the presence in phase space of a chaotic saddle, which is a fractal set. Orbits starting close to the chaotic saddle eventually escape in a highly unpredictable manner, and the escape times also depend on the disturbance and the initial conditions. In any case, one key feature of the method is that the control used is always smaller than the disturbance,  which is rather surprising and counterintuitive.

In the present work, we face a new objective, which can be viewed as the converse of the previous one. While in the previous case the goal was to keep the orbit in $Q$ for ever, now based on the same premises, our goal is to control the number of iterations necessary for the orbit to escape $Q$.

\begin{figure*}
	\centering
	\includegraphics[width=1\linewidth]{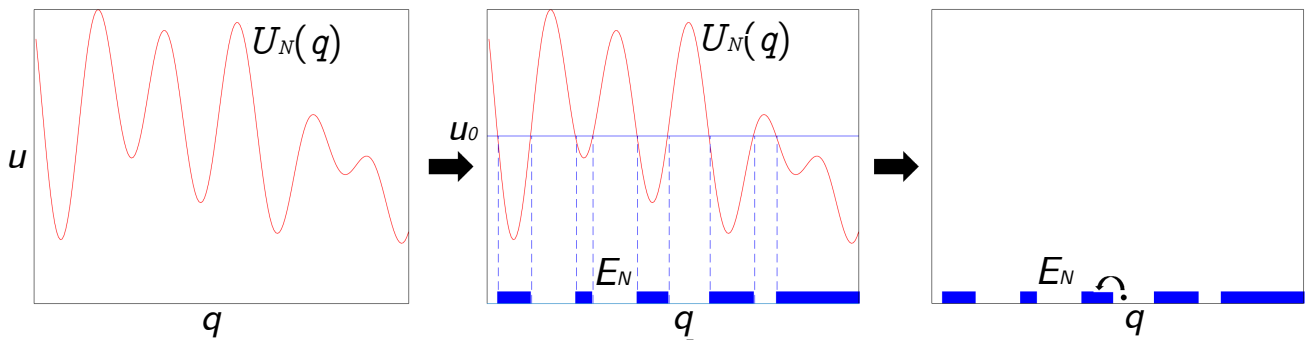}
	\caption{Steps that summarize the partial control procedure. (Left) Escape function $U_N(q)$. For a given initial point $q\in Q$, this function represents the minimum upper control bound $u_0$ necessary to expel an orbit out of the region $Q$ in $N$ iterations. (Center) Escape set $E_N$ (blue boxes) obtained as a result of the intersection from the escape function with the control bound $u_0$ (horizontal blue line). The points in $E_N$ will escape from $Q$ by using controls $|u_n| \leq u_0$. (Right) A point outside the escape set is controlled to insert it in the set. By doing so, each iteration guarantees that the orbit will escape in $N$ iterations or less.}
	\label{Flujo}
\end{figure*}

We have analyzed three different scenarios to apply the new strategy. The first one corresponds to the case where the controller wants to force the escape of the orbits from $Q$ in $N$ or less iterations of the map. In the second case, we consider the situation when we want the orbit to escape in exactly $N$ iterations, where this stronger constraint would necessarily imply a higher amount of control. Finally, the third situation we consider is somehow different. We have a map defined in two different regions in phase space, so that for a value of the parameter two chaotic attractors coexist independently in each region, and after a certain parameter value both attractors merge into a single global chaotic attractor. The idea here is to apply our control technique to fix the precise number of iterations of the orbit to stay in each region. As a result, we will get a chaotic orbit that periodically oscillates between the regions. Obviously this scheme could be generalized to a larger number of regions.

The paper is organized as follows. In Sect. 2, we introduce the partial control method and explain how to use it to control the escape from the chaotic region. Furthermore, we describe the escape functions and escape sets, adapted from previous work \cite{SafetyFun} that have been used for our objectives. In Sect. 3, we present two different specific cases where we apply the algorithm to fix the number of iterations for the orbits to escape. In Sect. 4, we address the case where the algorithm is applied for the goal of alternating the orbit between two regions in a predictable manner. Finally, the main conclusions are provided in the last section.

\section{The partial control method}

We present here the partial control method \cite{sasayo12,casasa18} that is applied on maps in the following manner
\begin{equation}
\begin{split}
&q_{n+1}=f(q_n)+\xi_n+u_n \\
&|\xi_n|\leq \xi_0 \\
&|u_n|\leq u_0 < \xi_0,
\end{split}
\end{equation}
where the map acts on values of $q\in Q$, $Q$ is a region in phase space, $\xi_n$ represents a bounded disturbance affecting the map at each iteration and $u_n$ is the applied control at each iteration, which is also bounded and importantly, smaller than the disturbance. 

\sloppy
Our goal here is to perturb the orbit starting in $Q$ by applying a sequence of controls $(u_1, u_2,.. ,u_N)$ in order to push the orbit out of $Q$ in $N$ iterations. Needless to say, we can achieve this objective by using different sequences, and the approach of partial control is to find the strategy that minimizes the upper bound of that sequence, that is, the $\min\big(\max\, (\,|u_1|, |u_2|,.. ,|u_N|\,)\,\big)=U_N(q)$, where $q$ is the initial point of the orbit and $U_N$ is the \textit{escape function}. Once the escape function is computed, we can choose an upper control bound value $u_0$ and select the set of points $q\in Q$ that satisfy $U_N(q) \leq u_0$. We name this set, \textit{the escape set} $E_N$. Any orbit starting in this set can be expelled from $Q$ by using a sequence of $N$ controls $u_n$ $ (n \leq N)$ with magnitude equal or smaller than $u_0$. The notions of escape functions and escape sets have been adapted from the safe functions and safe sets defined in \cite{SafetyFun}. The steps to apply this control technique are summarized as follows

\begin{enumerate}
\item  Choose the phase space region $Q$ where the control method will be applied. We assume that we know the map and the upper disturbance bound $\xi_0$ affecting it.
\item Compute the escape functions $U_k$ with $k=1:N$ in the region $Q$. Remind that $N$ is the number of iterations we need to expel the orbit out of $Q$.  
\item  Set the value $u_0$ and for every escape function $U_k$, compute the corresponding escape set $E_k$.
\item  For every iteration of the map, we choose the appropriate $u_n$, $|u_n|\leq u_0$, to bring the orbit to the  escape set $E_k$. Thus, we use the control so that the first iteration of the map brings the orbit to the escape set $E_N$, the second iteration of the map to the escape set $E_{N-1}$ and so on, until the orbit escapes out of the region $Q$.
\end{enumerate}

The steps $2$, $3$ and $4$ of this procedure are illustrated in Fig.~\ref{Flujo}. The second step is the most computationally expensive, where the escape functions $U_k$ are computed by using an algorithm developed in \cite{SafetyFun}. The algorithm is based on the observation that this kind of control problems can be solved backwards, starting from the last iteration. We will see that it is straightforward to compute the first function $U_1$ and through an iterative procedure obtain the rest of the escape functions $U_k$, since $U_{k+1}=f(U_k)$.

\begin{figure}[H]
	\centering
	\includegraphics[width=1\linewidth]{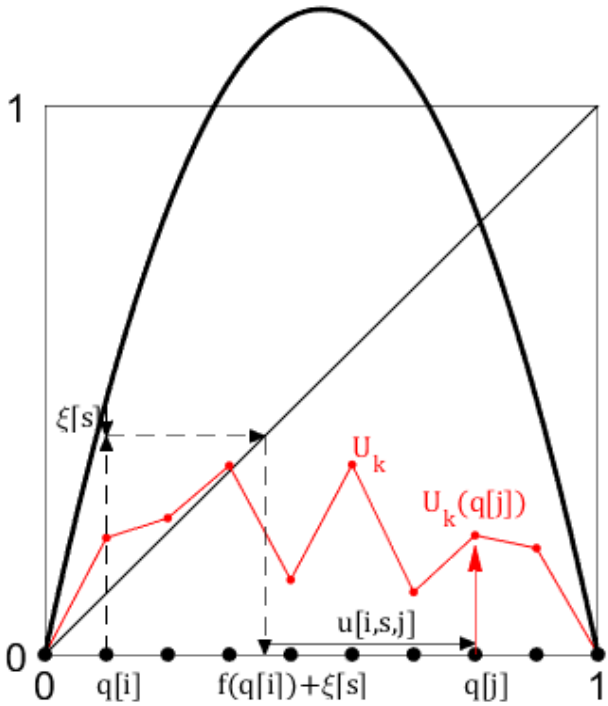}
	\caption{Scheme of the escape function (red) and the terminology used in the control procedure. Here the logistic map with $\mu=4.7$ is plotted. It is clear that orbits starting in points of $[0,1]$ will escape the interval $[0,1]$ after a few iterations. The dashed arrows are the mapping of an initial point $q[i]$ for one of the possible values of the perturbation $\xi[s]$. The horizontal black arrow represents the applied control that corresponds to the distance from the mapped point $f(q[i])+\xi[s]$ to the arrival point $q[j]=f(q[i])+\xi[s]+u[i,s,j]$. The indices $i$, $j$ refer to points in the discretization of the interval $[0,1]$, while $s$ refers to values of the discretization of $\xi_n$.}
	\label{EsquemaFun}
\end{figure}
Since we are doing numerical simulations, we must use a grid on $Q$ so that the map becomes
  \begin{equation}
  q_{n+1}=f(q_n)+\xi_n +u_n \;\;\; \rightarrow \;\;\; q[j]=f(q[i])+\xi[s]+u[i,s,j],
  \end{equation}
where $i=1:M$ denotes the number of the grid points in $Q$. The index $s=1:W$ corresponds with the number of possible disturbances ranging from $-\xi_0$ to $\xi_0$.  The index $j=1:M$ denotes the arrival point $q[j]=f(q[i])+\xi[s]+u[i,s,j]$. The term $u[i,s,j]$ denotes the control applied to the point $f(q[i],\xi[s])$ to put it in the arrival point $q[j]$. All these terms are illustrated in Fig.~\ref{EsquemaFun} for clarity.

\section{Escaping from the chaotic region}

We will use the well-known logistic map $x_{n+1}=f(x_n)=\mu x_n(1-x_n)$  as an example to illustrate the application of the algorithm described earlier. 

\begin{figure}
	\centering
	\includegraphics[width=1\linewidth]{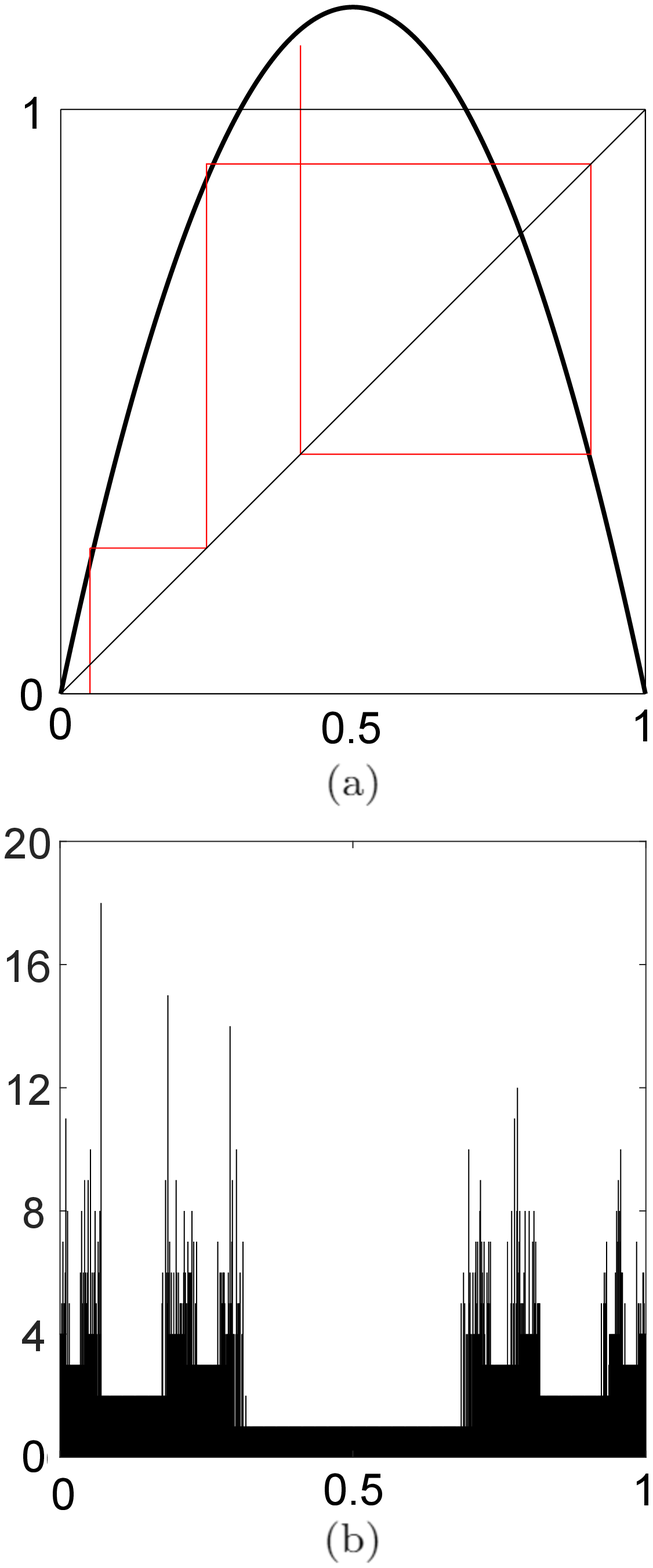}
	\caption{(a) Logistic map with $ \mu=4.7 $ and one orbit that eventually escapes from the region $Q=[0,1]$. (b) Escape time for the noise-affected logistic map. It clearly presents a fractal structure, so the time an orbit will escape is unpredictable. }
	\label{LogisticA}
\end{figure}

When we consider values of $\mu > 4$, the logistic map presents transient chaos in the region $Q=[0,1]$.  On Fig.~\ref{LogisticA}(a), the logistic map and one escaping orbit are represented. The number of iterations for which an orbit stays in $Q$ without control depends on the initial condition and the sequence of disturbances affecting it. As a consequence the lifetime is highly unpredictable as shown in Fig.~\ref{LogisticA}(b).

The goal of the algorithm based on the partial control method used in this work is to fix the number $N$ of the iterations after which the orbit escapes from $Q$. Three different cases will be explored. In each case, we will show the algorithm to compute the corresponding escape functions with some examples.

\subsection{Case A: Escape in N or less iterations}

This situation is justified in the case when we want an orbit to leave the chaotic region as quickly as possible. To do that, we choose the value $N$ and design the control algorithm so that the controlled orbit will abandon $Q$ in $N$ or less iterations. The lesser the value of $N$, the quicker the orbits will escape, though at a high price of control. 

As already commented, the escape functions $U_k$ are calculated from the first escape function $U_1$ through a recursive algorithm, since $U_{k+1}=f(U_k)$. Then, the value $U_1(q[i])$ is defined to be the minimum control necessary to escape in the next iteration. This control corresponds to the distance from the mapped point $f(q[i])+\xi[s]$ to the nearest end points of $Q$ (In this case $0$ or $1$). However, since there is a control $u[i,s,j]$ associated to each different value of the disturbance $\xi[s]$, we need to choose the maximum control among them all. Proceeding similarly for every initial condition $q[i]$, we obtain the escape function $U_1$. This function represents the minimum control bound necessary to escape in one iteration.

The next function $U_2$ will correspond to the minimum control bound necessary to force the escape of the orbit within $2$ iterations. There are two possibilities here, given an initial condition $q[i]$, the image $f(q[i])+\xi[s]$ can fall directly outside $Q$, in which case no control is needed, or can fall in a point $q[j] \in Q$. In the latter case, we need to compute the suitable control that minimizes the maximum between the values $u[i,s,j]$ (the control applied in this iteration of the map) and $U_1(q[j])$ (the maximum control that we will apply in the next iteration of the map), that is, the $\min\limits_{j}(\max(u[i,s,j],U_1[j]))$. Again, to take into account all possible disturbances $\xi[s]$ we must choose the maximum among all corresponding controls, that is $U_2(q[i])=\max\limits_{s}(\min\limits_{j}(\max(u[i,s,j],U_1(q[j]))))$.
 
The procedure to compute $U_3$ is similar, where now  $U_3(q[i])=\max\limits_{s}(\min\limits_{j}(\max(u[i,s,j],U_2(q[j]))))$. The algorithm is repeated until $U_N$ is obtained. The values $U_N(q[i])$ of this function represent the minimum control bound necessary to force the orbit starting in $q[i]$ to escape from $Q$ within $N$ or less iterations. 

In order to define the algorithm to compute the escape functions, we define $u_{out}[i,s]$ as the control applied to the image $f(q[i], \xi[s])$ to move it outside the region $Q$. Then, $U_1$ is calculated as follows. 

\begin{equation}
U_1(q[i])=\max_s(u_{out}[i,s])
\end{equation}
If $q[i] + \xi[s]$ is beyond $Q$, then $u_{out}[i,s]=0$. Then, given $U_1$ as the seed function, we can calculate the next escape functions with the following recursive algorithm

\begin{equation}\small
U^{in}_{k+1}(q[i])=\max_s\Big(\min_j\big(\max(u[i,s,j],U_k(q[j]))\big)\Big)
\label{Uaux}
\end{equation}
\begin{equation}\small
U_{k+1}(q[i])=\min\bigg(U_1(q[i]),U^{in}_{k+1}(q[i])\bigg),
\end{equation}
where the intermediate function $U^{in}_{k+1}$ was introduced to allow the orbit to escape from $Q$ before $N$ iterations. Notice that $U^{in}_{k+1}$ only takes into account images inside $Q$ to control the orbit, while  $U_1$ only takes into account images outside $Q$. Between these two possibilities, the one that minimizes the control will be chosen. By doing so, the orbit can be expelled in any of the $k \leq N$ iterations.

\begin{figure}
	\centering
	\includegraphics[width=1\linewidth]{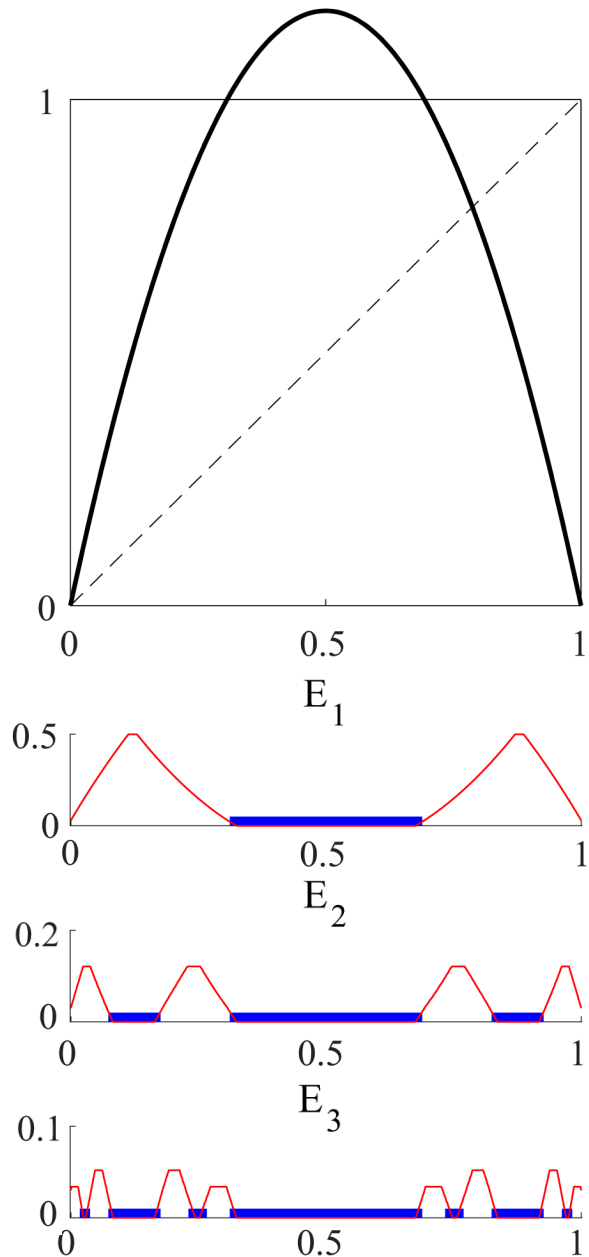}
	\caption{Case A. On top, the logistic map with $ \mu=4.7 $. On the bottom, the escape functions $U_k$ are computed for $N=3$ and shown in red. The escape sets $E_k$ are shown in blue. We have used $ \xi_0=0.030 $ and $u_0=0.022 $. Controlled orbits starting in $E_3$ will escape from $Q=[0,1]$ in $3$ or less iterations by applying a control $|u_n|\leq u_0$ at each iteration. The scales used in the vertical axis are different for a better visualization.}
	\label{EscapeSetA}
\end{figure}

As an example, we have chosen $N=3$ so that controlled orbits will escape from $Q$ in $3$ or less iterations. The upper disturbance bound affecting the logistic map was set to $\xi_0= 0.030$. The corresponding escape functions $U_k$ are shown in Fig.~\ref{EscapeSetA}, which as it can be observed take zero values in some intervals. This means that points in these intervals where $U_k=0$ will escape from $Q$ within $k$ iterations for any $\xi[s]$, without applying any control. It can be also observed that as the index $k$ increases, the escape function decreases since the orbit has more iterations to escape.

\begin{figure}
	\centering
	\includegraphics[width=1\linewidth]{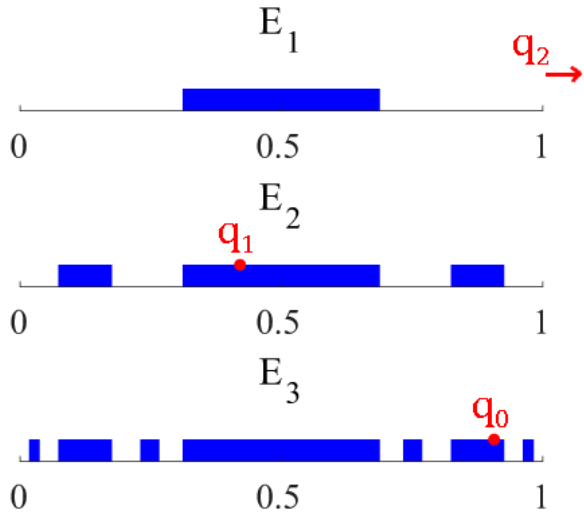}
	\caption{Case A. Example of a controlled orbit. An initial condition in $E_3$ is controlled to $E_2$, then  $E_1$, and finally it leaves the chaotic region $Q=[0,1]$ under a suitable control. In general, orbits can escape in $3$ or less iterations depending on the initial condition and the disturbance $\xi_n$. At each iteration, the applied control $u_n$ is the minimum between the nearest escape set and the end points of $Q$. We have fixed here $\mu=4.7$, $\xi_0=0.030$ and $u_0=0.022$}
	\label{TrayA}
\end{figure}

Once we compute the escape functions, we have to select the control value $u_0$ to compute the corresponding escape sets $E_k$. No tall $u_0$ values are allowed. Escape sets only exist for values $u_0 \geq min(U_N)$. The bigger $u_0$, the bigger the escape sets. These sets consist of points $q[i]$ satisfying the condition $ U_k(q[i]) \leq u_0 $. The escape sets $E_k$ for $u_0 = 0.022$ are represented in Fig.~\ref{EscapeSetA}. For this case, the set $E_3$ represents the set of points $q[i]$ that can escape from $Q$ within $3$ iterations of the map, by applying a control $u_n \leq u_0$ at each iteration. In the first iteration, the control $u_1$ will be applied to put the orbit in the closest point of $E_2$ or outside $Q$ if possible. In the second iteration, the control $u_2$ will be applied to put the orbit in the closest point of $E_1$ or outside $Q$ if possible. In the third and last iteration the control $u_3$ will be applied to put the orbit outside $Q$. This is illustrated in Fig.~\ref{TrayA}.

\begin{figure}
	\centering
	\includegraphics[width=1\linewidth]{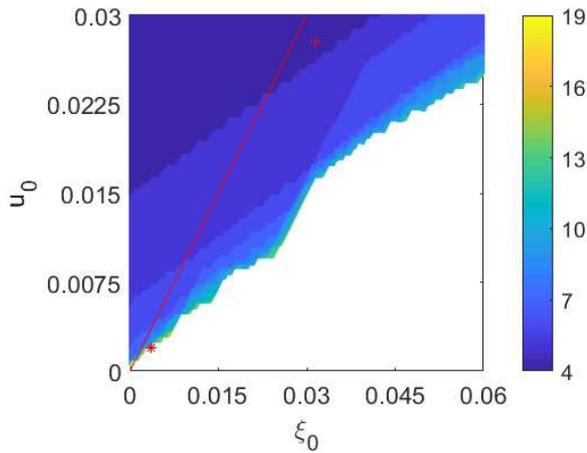}
	\caption{This color plot shows the minimum number of iterations $N$ needed for every orbit to escape from $Q=[0,1]$ for any given couple of parameters ($\xi_0$, $u_0$), while the white points correspond to a situation where some initial conditions cannot escape. The red line represents the ratio $u_0=\xi_0$. As an illustration we have marked two points in red. The (+) point corresponds to $N=4$, while the (*) point to $N=19$.}
	\label{Ncolor}
\end{figure}

We have built in Fig.~\ref{Ncolor} a colormap plot showing the minimum number of iterations $N$ needed for every orbit to escape from $Q=[0,1]$ in a ($\xi_0$, $u_0$) parameter plane. In this figure the white points correspond to the case where some initial conditions cannot escape from $Q$ using these values of ($\xi_0$, $u_0$). A red line is plotted representing the ratio $u_0=\xi_0$, so that values below this line are the ones we are interested due to the partial control method. Furthermore, we have chosen two points marked in red. The (+) point corresponds to $N=4$, while the (*) point to $N=19$.

\subsection{Case B: Escape exactly in N iterations}

Here we analyze the case when the number of iterations for an orbit to escape from $Q$ is exactly $N$. A relevant observation here is that we need to control the orbit inside $Q$ for $N-1$ iterations, since it will escape precisely at the iteration $N$. The algorithm to obtain the escape functions is now simpler, because we do not need to consider the possibility that the orbit abandons $Q$ before the iteration $N$, and is described next

\begin{equation}
U_1(q[i])=\max_s(u_{out}[i,s])
\end{equation}

\begin{equation}\small
U_{k}(q[i])=\max_s\Big(\min_j\big(\max(u[i,s,j],U_k(q[j]))\big)\Big).
\end{equation}

\begin{figure}
	\centering
	\includegraphics[width=1\linewidth]{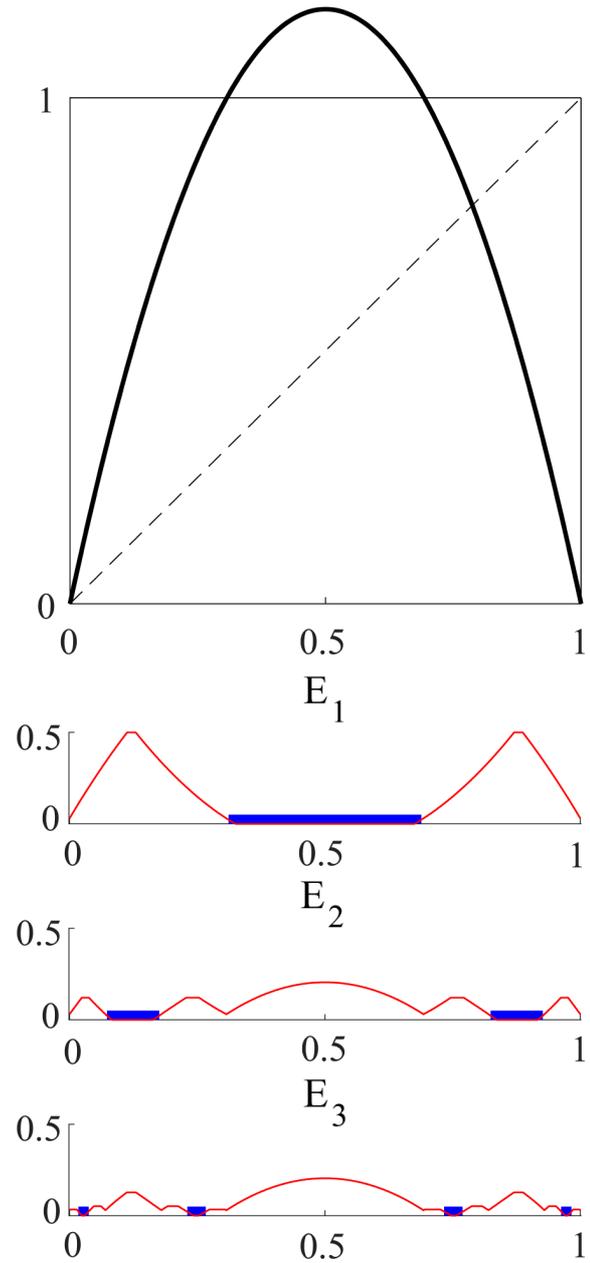}
	\caption{Case B. On top, the logistic map with $ \mu=4.7 $. On the bottom, the escape functions $U_k$ are computed for $N=3$ and shown in red.  The escape sets $E_k$ are shown in blue. We have used $ \xi_0=0.030 $ and $u_0=0.022 $. Orbits starting in $E_3$ will escape from $Q=[0,1]$ in exactly $3$ iterations with a control $|u_n|\leq u_0$ at each iteration. }
	\label{EscapeSetB}
\end{figure}

To show how to control orbits to escape from $Q$ in exactly $3$ iterations, we use the logistic map with $ \mu=4.7 $, and $ \xi_0=0.030 $. The escape functions are represented in Fig.~\ref{EscapeSetB}. As it can be observed, the escape functions in case B have bigger values than in case A. The reason could be that the condition ``expel the orbit in exactly $N$ iterations" (case B) is stronger than the condition ``expelling the orbit in $N$ or less iterations" (case A). Therefore, we need bigger controls for the case B than the case A.

\begin{figure}
	\centering
	\includegraphics[width=1\linewidth]{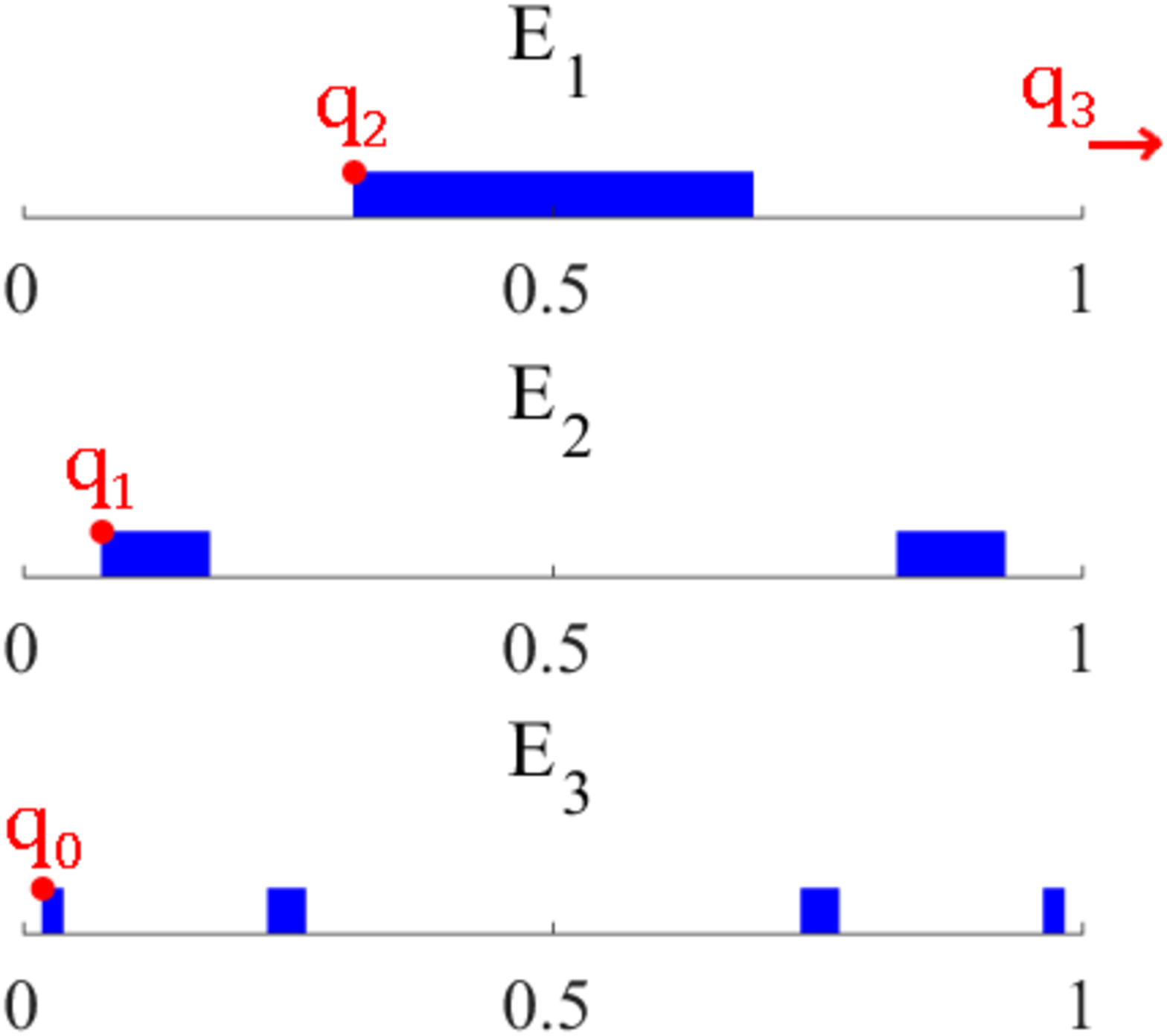}
	\caption{Case B. Example of a controlled orbit. An initial condition in $E_3$ is mapped to $E_2$, then  $E_1$, and finally it leaves the chaotic region $Q=[0,1]$ after a suitable control. All orbits starting in $E_3$ will escape $Q=[0,1]$ in exactly $3$ iterations.}
	\label{TrayB}
\end{figure}

To get the escape sets $E_k$, we have chosen the control bound $u_0=0.022$. These sets are represented in Fig.~\ref{EscapeSetB}. An scheme of how an orbit starting in $E_3$ escapes $Q=[0,1]$ in exactly $3$ iterations is shown in Fig.~\ref{TrayB}.

\section{Case C: Alternating chaotic regions}

In the previous section, we mentioned three different cases to explore the escaping of an orbit from $Q$ in $N$ iterations. We discussed earlier, cases A and B. Now, we focus our attention on the third case corresponding to multistable chaotic systems that merge into a larger chaotic attractor as a parameter is varied \cite{rempel2005intermittency,livorati2015crises,vaidyanathan2018new}. This process occurs mainly when the basin boundary of each attractor collide, so that an orbit moves chaotically back and forth from one region to the other \cite{jin2018stochastic}.

Our goal here is to use the control algorithm to force the orbit to stay in each region for a fixed number of iterations before moving to the other region. What we want here is to maintain a perpetual periodical motion bouncing back and forth between regions. The initial chaotic orbit stays for $N$ iterations in one region, then it moves to the second region where it stays for $N'$ iterations, and finally it comes back again to the first region. As we will see, the orbit will resemble a chaotic signal modulated by a periodic one.   

\begin{figure}
	\centering
	\includegraphics[width=1\linewidth]{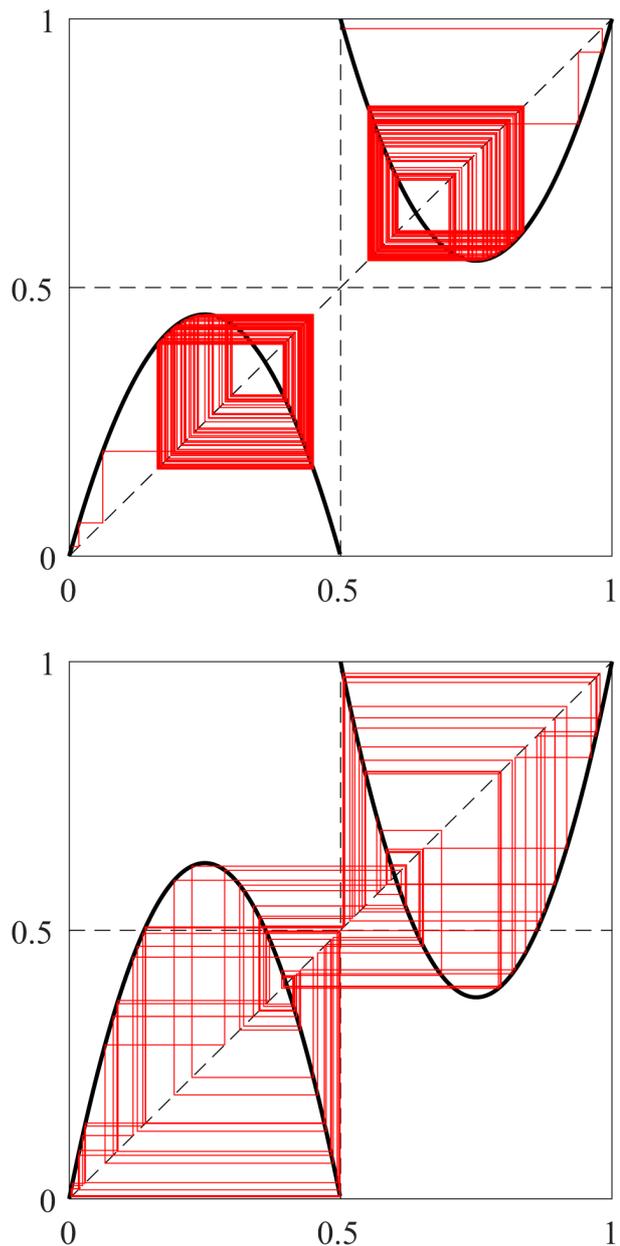}
	\caption{\textbf{Double parabola map}. The map is defined in a piecewise manner in two different regions $Q_l=[0,0.5)$ and $Q_r=[0.5,1]$. (Top) When $ \mu = 7.2 $, the orbits stay on one side of the map, so there are clearly two different attractors. (Bottom) When $ \mu = 10 $, a global chaotic attractor merges. Now, orbits starting in $Q_l$ have a transient chaotic behavior before escaping to the region $Q_r$ where after another chaotic transient the orbit comes back to $Q_l$.}
	\label{Alt_Map}
\end{figure}

We have constructed a map that we name \textit{the double parabola map} illustrated in Fig.~\ref{Alt_Map} as an example of a simple map exhibiting the behavior described before, 

\begin{equation}\small
x_{n+1} = \left\{ \begin{array}{ll}
-\mu (x_n^2 + \dfrac{1}{2}x_n)  & \mbox{si $x<0.5$,} \\
1 + \mu (x_n^2 -\dfrac{3}{2}x_n + \dfrac{1}{2})  & \mbox{si $x_n\geq 0.5$,} 
\end{array}
\right.
\end{equation}

This map is defined by a convex parabola at the region of the left $Q_l=[0,0.5)$ and a concave parabola at the right region $Q_r=[0.5,1]$. Both parabolas depend on the same parameter $\mu$. 

To show an example of the application of our logarithm to this case, we will focus on the behavior of the map for $\mu=10$, where we have a chaotic attractor that expands to all the interval $[0,1]$. Thus, orbits starting in the region $Q_l$ behave chaotically to eventually escaping to the region $Q_r$, and vice versa as shown in the Fig.~\ref{Alt_Map}. Even though the map is well defined to map points of $[0,1]$ into itself, however this could not be so in presence of a disturbance. In particular, for points close to the end of interval, what it should be considered in the control scheme.

The main purpose now is to apply a control so that the transition from the two regions would be predictable. In other words, we want to control how many iterations a given orbit stays on each region. As a consequence, our goal here is to keep the orbit $N_l$ iterations on region $Q_l$ and $N_r$ iterations on region $Q_r$, where the values $N_l$ and $N_r$ are previously chosen by the controller. As a result, we will get a chaotic orbit that periodically oscillates between regions $Q_l$ and $Q_r$. 

To compute the escape functions, we follow a similar methodology as the one used on case B from the previous section, though we must impose in the algorithm the periodic condition.  Here, we will have $N_l$ escape functions for orbits in $Q_l$ and $N_r$ escape functions for orbits in $Q_r$. From now on, we will denote as $U^l$ and $E^l$ the escape functions and escape sets of the left region $Q_l$. Similarly $U^r$ and $E^r$ are defined in the right region $Q_r$. 

Next, we briefly describe the algorithm to compute the escape functions. First, we start by computing the function

\begin{equation}
U^l_1(q[i])=\max\limits_{s}\Big(\min_j\big(u_{in}^l[i,s,j]\big)\Big).
\end{equation}

Taking this function as a seed, we can compute the rest of the escape functions with the following algorithm

\SetKwBlock{DummyBlock}{}{}

\begin{algorithm}
\SetAlgoLined
\mbox{}\\
\mbox{}\\

Loop until the functions $U^{l}_{k}$ and $U^{l}_{k}$ converges
\DummyBlock{\SetAlgoLined
\For{$k=1:N_l-1$} {
    $U^{l}_{k+1}(q[i])=\max\limits_s\Big(\min\limits_j\big(\max(u_{in}^l[i,s,j],U^l_k(q[j]))\big)\Big)$	
    }
$U^{r}_{1}(q[i])=\max\limits_s\Big(\min\limits_j\big(\max(u_{out}^r[i,s,j],U^l_{N_l}(q[j]))\big)\Big)$

\For{$k=1:N_r-1$} {
$U^{r}_{k+1}(q[i])=\max\limits_s\Big(\min\limits_j\big(\max(u_{in}^r[i,s,j],U^r_k(q[j]))\big)\Big)$
}
$U^{l}_{1}(q[i])=\max\limits_s\Big(\min\limits_j\big(\max(u_{out}^l[i,s,j],U^r_{N_r}(q[j]))\big)\Big)$
}
\mbox{}\\
\mbox{}\\
\end{algorithm}

\noindent
where $ u_{in}^{l} $ denotes the control applied to remain in $Q_l$ and  $ u_{in}^{r} $ is the one to remain in $Q_r$. On the other hand, $ u_{out}^{r} $ denotes the control needed to migrate from $Q_r$ to $Q_l$ and conversely $ u_{out}^{l} $ is the control to migrate from $Q_l$ to $Q_r$.

\begin{figure}
	\centering
	\includegraphics[width=1\linewidth]{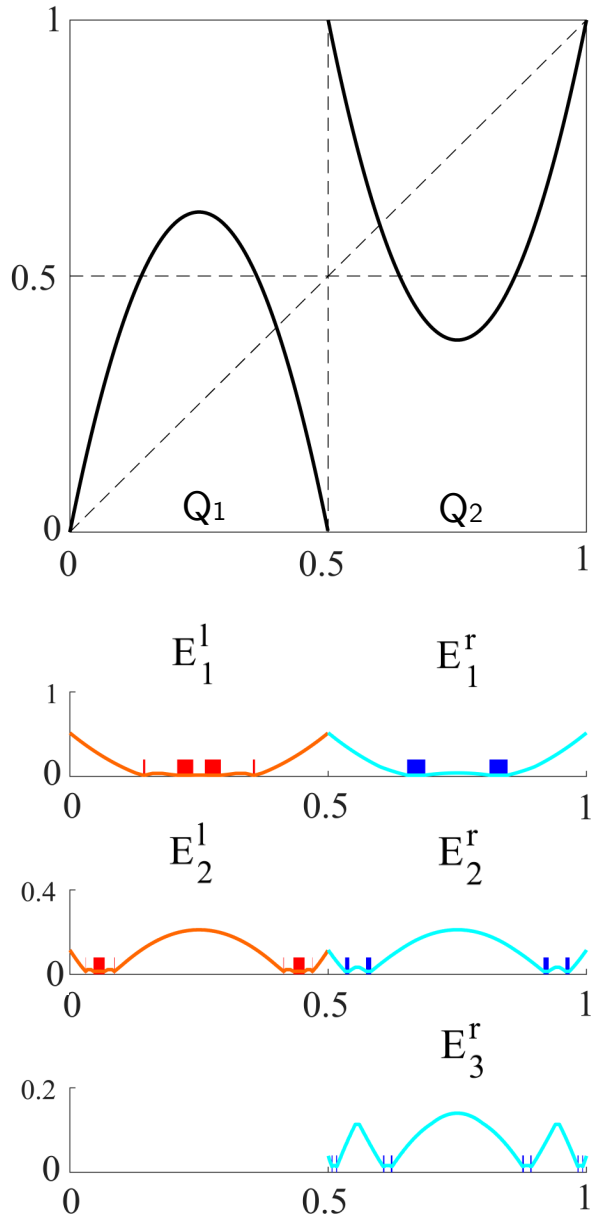}
	\caption{On the top, the double parabola map with $ \mu=10 $. On the bottom, the escape functions computed for $ N_l=2 $ (orange) and  $ N_r=3 $ (cyan). Furthermore, the escape sets $E^{r}$ appear as blue boxes and $E^{l}$ as red boxes. We have used $ \xi_0=0.015 $ as the disturbance bound and $ u_0=0.014 $ as the control bound, where this one corresponds to the minimum value of the escape functions.} 
	\label{AltSet}
\end{figure}

\begin{figure}
	\centering
	\includegraphics[width=1\linewidth]{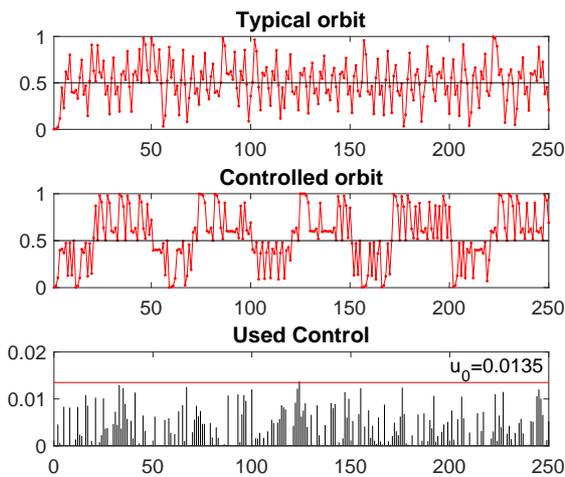}
	\caption{On the top, a typical orbit for the double parabola map ($ \mu=10 $). In the middle, we show an orbit affected by a disturbance $ \xi_0=0.0150 $ where we have used an upper control bound $u_0=0.0135$ to keep the orbit  $N^l=20$ iterations in $Q_l=[0,0.5]$ and $N^r=30$ iterations in $Q_r=[0.5,1]$. At the bottom, the absolute value of the control used during the $250$ iterations is represented. Note that all the applied controls remain below the control bound $u_0=0.0135$ shown by the red line, which corresponds to the minimum value of the escape functions. }
	\label{AltTray}
\end{figure}

Now, we want to illustrate the computation of a simple example where we have chosen  $N_l=2$, $N_r=3$ and the upper disturbance bound $ \xi_0=0.015 $. The computed escape functions  $U^{l}$ and  $U^{r}$ are shown in Fig.~\ref{AltSet}. For the computation of the escape sets $E_k$, we choose the control bound $u_0=0.014$, which corresponds to the minimum value of the escape functions, and compute the points satisfying the conditions $ U^{l}_k\leq u_0 $ and $ U^{r}_k\leq u_0 $ so that we obtain $E^{l}_k$ and $E^{r}_k$, respectively, as shown in Fig.~\ref{AltSet}. To control the orbit, we need to choose an initial condition in one of the escape sets $E^{l}_k$ or $E^{r}_k$, and then a suitable control is applied. For example, if we start with a point in $E^{r}_3$, we apply the control in the next iterations to put the orbit in $E^{r}_2 \rightarrow E^{r}_1 \rightarrow E^{l}_2 \rightarrow E^{l}_1 \rightarrow E^{r}_3 \rightarrow E^{r}_2 \rightarrow ...(repeat)$. As a result, we obtain a chaotic motion modulated by a periodic one.

To make this behavior even more clear, we show another example considering now $N^l=20$, $N^r=30$, $ \xi_0=0.015 $ and control bound $u_0=0.0135$. This is shown in Fig.~\ref{AltTray}, on top it appears a typical orbit of the map and in the middle the controlled orbit remaining $20$ iterations in $Q_l=[0,0.5)$  and $30$ 
iterations in $Q_r=[0.5,1)$. On the bottom, the values of the applied control are shown, which are all below $u_0=0.0135$.

\section{Conclusions}

We have developed a new control algorithm based on the partial control method aiming at keeping an orbit on a certain region $Q$ for a given number of iterations $N$ with a minimum control. For that purpose, we have adapted known tools such as safe functions and safe sets to the new escape functions and the escape sets. Once the latter are computed, it is straightforward to control the orbit. 

We have considered three different possible scenarios. The first case, case A, where the main goal has been to force the escape of a given orbit of the phase space $Q$ in $N$ or less iterations. The second case, case B, where the goal has been to force the escape of the orbit to happen in exactly $N$ iterations. 

And the third scenario, case C, where we use the algorithm in a situation where a system has chaotic transitions between two regions, with the goal to control these transitions. As a consequence, this allows the controller to fix the number of iterations that the orbit will stay in every particular region, making it to have a periodic sequence of transitions. 

Even though we have used one-dimensional maps for simplicity, we believe our control method is valid for higher dimensions, and we hope that it can be applied to different problems modeled with maps.

\begin{acknowledgements}
This work was supported by the Spanish State Research Agency (AEI) and the European Regional Development Fund (ERDF, EU) under Project No.~PID2019-105554GB-I00.
\end{acknowledgements}

\section*{Compliance with ethical standards}
\textbf{Conflict of interest} The authors declare that they have no conflict of interest concerning the publication of this manuscript.

\end{document}